\documentclass[11pt]{amsart}
\usepackage{amsmath}
\usepackage{amssymb,latexsym}
\usepackage{amscd}
\usepackage{xspace}
\usepackage{verbatim}
\begin{document}
\numberwithin{equation}{section}
\title[Reduced limit for  semilinear problems with measures]{Reduced limit for  semilinear boundary value problems with measure data.}
\author{Mousomi Bhakta and Moshe Marcus}
\address{Department of Mathematics, Technion\\
 Haifa 32000, ISRAEL}
\date{\today}

\newcommand{\txt}[1]{\;\text{ #1 }\;}
\newcommand{\tbf}{\textbf}
\newcommand{\tit}{\textit}
\newcommand{\tsc}{\textsc}
\newcommand{\trm}{\textrm}
\newcommand{\mbf}{\mathbf}
\newcommand{\mrm}{\mathrm}
\newcommand{\bsym}{\boldsymbol}
\newcommand{\scs}{\scriptstyle}
\newcommand{\sss}{\scriptscriptstyle}
\newcommand{\txts}{\textstyle}
\newcommand{\dsps}{\displaystyle}
\newcommand{\fnz}{\footnotesize}
\newcommand{\scz}{\scriptsize}
\newcommand{\be}{\begin{equation}}
\newcommand{\bel}[1]{\begin{equation}\label{#1}}
\newcommand{\ee}{\end{equation}}
\newtheorem{subn}{\name}
\renewcommand{\thesubn}{}
\newcommand{\bsn}[1]{\def\name{#1$\!\!$}\begin{subn}}
\newcommand{\esn}{\end{subn}}
\newtheorem{sub}{\name}[section]
\newcommand{\dn}[1]{\def\name{#1}}   
\newcommand{\bs}{\begin{sub}}
\newcommand{\es}{\end{sub}}
\newcommand{\bsl}[1]{\begin{sub}\label{#1}}
\newcommand{\bth}[1]{\def\name{Theorem}\begin{sub}\label{t:#1}}
\newcommand{\blemma}[1]{\def\name{Lemma}\begin{sub}\label{l:#1}}
\newcommand{\bcor}[1]{\def\name{Corollary}\begin{sub}\label{c:#1}}
\newcommand{\bdef}[1]{\def\name{Definition}\begin{sub}\label{d:#1}}
\newcommand{\bprop}[1]{\def\name{Proposition}\begin{sub}\label{p:#1}}
\newcommand{\bnote}[1]{\def\name{\mdseries\scshape Notation}\begin{sub}\label{n:#1}}
\newcommand{\bproof}{\begin{proof}}
\newcommand{\eproof}{\end{proof}}
\newcommand{\bcom}{}
\newcommand{\req}{\eqref}
\newcommand{\rth}[1]{Theorem~\ref{t:#1}}
\newcommand{\rlemma}[1]{Lemma~\ref{l:#1}}
\newcommand{\rcor}[1]{Corollary~\ref{c:#1}}
\newcommand{\rdef}[1]{Definition~\ref{d:#1}}
\newcommand{\rprop}[1]{Proposition~\ref{p:#1}}
\newcommand{\rnote}[1]{Notation~\ref{n:#1}}
\newcommand{\BA}{\begin{array}}
\newcommand{\EA}{\end{array}}
\newcommand{\BAN}{\renewcommand{\arraystretch}{1.2}
\setlength{\arraycolsep}{2pt}\begin{array}}
\newcommand{\BAV}[2]{\renewcommand{\arraystretch}{#1}
\setlength{\arraycolsep}{#2}\begin{array}}
\newcommand{\BSA}{\begin{subarray}}
\newcommand{\ESA}{\end{subarray}}
\newcommand{\BAL}{\begin{aligned}}
\newcommand{\EAL}{\end{aligned}}
\newcommand{\BALG}{\begin{alignat}}
\newcommand{\EALG}{\end{alignat}}
\newcommand{\BALGN}{\begin{alignat*}}
\newcommand{\EALGN}{\end{alignat*}}
\newcommand{\note}[1]{\noindent\textit{#1.}\hspace{2mm}}
\newcommand{\Remark}{\note{Remark}}
\newcommand{\forevery}{\quad \forall}
\newcommand{\1}{\\[1mm]}
\newcommand{\2}{\\[2mm]}
\newcommand{\3}{\\[3mm]}
\newcommand{\set}[1]{\{#1\}}
\def\capa{\mathit{cap}}
\newcommand{\st}[1]{{\rm (#1)}}
\newcommand{\lra}{\longrightarrow}
\newcommand{\lla}{\longleftarrow}
\newcommand{\llra}{\longleftrightarrow}
\newcommand{\Lra}{\Longrightarrow}
\newcommand{\Lla}{\Longleftarrow}
\newcommand{\Llra}{\Longleftrightarrow}
\newcommand{\warrow}{\rightharpoonup}
\def\dar{\downarrow}
\def\uar{\uparrow}
\newcommand{\paran}[1]{\left (#1 \right )}
\newcommand{\sqrbr}[1]{\left [#1 \right ]}
\newcommand{\curlybr}[1]{\left \{#1 \right \}}
\newcommand{\absol}[1]{\left |#1\right |}
\newcommand{\norm}[1]{\left \|#1\right \|}
\newcommand{\angbr}[1]{\left< #1\right>}
\newcommand{\paranb}[1]{\big (#1 \big )}
\newcommand{\sqrbrb}[1]{\big [#1 \big ]}
\newcommand{\curlybrb}[1]{\big \{#1 \big \}}
\newcommand{\absb}[1]{\big |#1\big |}
\newcommand{\normb}[1]{\big \|#1\big \|}
\newcommand{\angbrb}[1]{\big\langle #1 \big \rangle}
\newcommand{\thkl}{\rule[-.5mm]{.3mm}{3mm}}
\newcommand{\thknorm}[1]{\thkl #1 \thkl\,}
\newcommand{\trinorm}[1]{|\!|\!| #1 |\!|\!|\,}
\newcommand{\vstrut}[1]{\rule{0mm}{#1}}
\newcommand{\rec}[1]{\frac{1}{#1}}
\newcommand{\opname}[1]{\mathrm{#1}\,}
\newcommand{\supp}{\opname{supp}}
\newcommand{\dist}{\opname{dist}}
\newcommand{\sign}{\opname{sign}}
\newcommand{\diam}{\opname{diam}}
\newcommand{\q}{\quad}
\newcommand{\qq}{\qquad}
\newcommand{\hsp}[1]{\hspace{#1mm}}
\newcommand{\vsp}[1]{\vspace{#1mm}}
\newcommand{\prt}{\partial}
\newcommand{\sms}{\setminus}
\newcommand{\ems}{\emptyset}
\newcommand{\ti}{\times}
\newcommand{\pr}{^\prime}
\newcommand{\ppr}{^{\prime\prime}}
\newcommand{\tl}{\tilde}
\newcommand{\wtl}{\widetilde}
\newcommand{\sbs}{\subset}
\newcommand{\sbeq}{\subseteq}
\newcommand{\nind}{\noindent}
\newcommand{\ovl}{\overline}
\newcommand{\unl}{\underline}
\newcommand{\nin}{\not\in}
\newcommand{\pfrac}[2]{\genfrac{(}{)}{}{}{#1}{#2}}
\newcommand{\tin}{\to\infty}
\newcommand{\ind}[1]{_{_{#1}}}
\newcommand{\chr}[1]{\mathbf{1}\ind{#1}}
\newcommand{\rest}[1]{\big |\ind{#1}}
\newcommand{\Sol}[2]{\mathrm{Sol}\ind{#2}{#1}}
\newcommand{\wkc}{weak convergence\xspace}
\newcommand{\wrto}{with respect to\xspace}
\newcommand{\cons}{consequence\xspace}
\newcommand{\consy}{consequently\xspace}
\newcommand{\Consy}{Consequently\xspace}
\newcommand{\Essy}{Essentially\xspace}
\newcommand{\essy}{essentially\xspace}
\newcommand{\sth}{such that\xspace}
\newcommand{\ngh}{neighborhood\xspace}
\newcommand{\seq}{sequence\xspace}
\newcommand{\sseq}{subsequence\xspace}
\newcommand{\ifif}{if and only if\xspace}
\newcommand{\suff}{sufficiently\xspace}
\newcommand{\abc}{absolutely continuous\xspace}
\newcommand{\subs}{sub-solution\xspace}
\newcommand{\supers}{super-solution\xspace}
\newcommand{\Wlg}{Without loss of generality\xspace}
\newcommand{\wlg}{without loss of generality\xspace}
\newcommand{\locun}{locally uniformly\xspace}
\newcommand{\bvp}{boundary value problem\xspace}
\newcommand{\bdw}{\partial\Gw}
\newcommand{\tr}{\mathrm{tr}\,}
\newcommand{\Tr}{\mathrm{Tr}\,}
\newcommand{\gsmod}{$\gs$-moderate\xspace}
\newcommand{\ofrown}{\overset{\frown}}
\def\RN{\BBR^N}
\def\loc{\ind{\rm loc}}
\def\bcom{}
\def\ga{\alpha}     \def\gb{\beta}       \def\gg{\gamma}
\def\gc{\chi}       \def\gd{\delta}      \def\ge{\epsilon}
\def\gth{\theta}                         \def\vge{\varepsilon}
\def\gf{\phi}       \def\vgf{\varphi}    \def\gh{\eta}
\def\gi{\iota}      \def\gk{\kappa}      \def\gl{\lambda}
\def\gm{\mu}        \def\gn{\nu}         \def\gp{\pi}
\def\vgp{\varpi}    \def\gr{\rho}        \def\vgr{\varrho}
\def\gs{\sigma}     \def\vgs{\varsigma}  \def\gt{\tau}
\def\gu{\upsilon}   \def\gv{\vartheta}   \def\gw{\omega}
\def\gx{\xi}        \def\gy{\psi}        \def\gz{\zeta}
\def\Gg{\Gamma}     \def\Gd{\Delta}      \def\Gf{\Phi}
\def\Gth{\Theta}
\def\Gl{\Lambda}    \def\Gs{\Sigma}      \def\Gp{\Pi}
\def\Gw{\Omega}     \def\Gx{\Xi}         \def\Gy{\Psi}

\def\CS{{\mathcal S}}   \def\CM{{\mathcal M}}   \def\CN{{\mathcal N}}
\def\CR{{\mathcal R}}   \def\CO{{\mathcal O}}   \def\CP{{\mathcal P}}
\def\CA{{\mathcal A}}   \def\CB{{\mathcal B}}   \def\CC{{\mathcal C}}
\def\CD{{\mathcal D}}   \def\CE{{\mathcal E}}   \def\CF{{\mathcal F}}
\def\CG{{\mathcal G}}   \def\CH{{\mathcal H}}   \def\CI{{\mathcal I}}
\def\CJ{{\mathcal J}}   \def\CK{{\mathcal K}}   \def\CL{{\mathcal L}}
\def\CT{{\mathcal T}}   \def\CU{{\mathcal U}}   \def\CV{{\mathcal V}}
\def\CZ{{\mathcal Z}}   \def\CX{{\mathcal X}}   \def\CY{{\mathcal Y}}
\def\CW{{\mathcal W}}
\def\BBA {\mathbb A}   \def\BBb {\mathbb B}    \def\BBC {\mathbb C}
\def\BBD {\mathbb D}   \def\BBE {\mathbb E}    \def\BBF {\mathbb F}
\def\BBG {\mathbb G}   \def\BBH {\mathbb H}    \def\BBI {\mathbb I}
\def\BBJ {\mathbb J}   \def\BBK {\mathbb K}    \def\BBL {\mathbb L}
\def\BBM {\mathbb M}   \def\BBN {\mathbb N}    \def\BBO {\mathbb O}
\def\BBP {\mathbb P}   \def\BBR {\mathbb R}    \def\BBS {\mathbb S}
\def\BBT {\mathbb T}   \def\BBU {\mathbb U}    \def\BBV {\mathbb V}
\def\BBW {\mathbb W}   \def\BBX {\mathbb X}    \def\BBY {\mathbb Y}
\def\BBZ {\mathbb Z}

\def\GTA {\mathfrak A}   \def\GTB {\mathfrak B}    \def\GTC {\mathfrak C}
\def\GTD {\mathfrak D}   \def\GTE {\mathfrak E}    \def\GTF {\mathfrak F}
\def\GTG {\mathfrak G}   \def\GTH {\mathfrak H}    \def\GTI {\mathfrak I}
\def\GTJ {\mathfrak J}   \def\GTK {\mathfrak K}    \def\GTL {\mathfrak L}
\def\GTM {\mathfrak M}   \def\GTN {\mathfrak N}    \def\GTO {\mathfrak O}
\def\GTP {\mathfrak P}   \def\GTR {\mathfrak R}    \def\GTS {\mathfrak S}
\def\GTT {\mathfrak T}   \def\GTU {\mathfrak U}    \def\GTV {\mathfrak V}
\def\GTW {\mathfrak W}   \def\GTX {\mathfrak X}    \def\GTY {\mathfrak Y}
\def\GTZ {\mathfrak Z}   \def\GTQ {\mathfrak Q}
\font\Sym= msam10
\def\SYM#1{\hbox{\Sym #1}}
\def\bmn{\mathbf{n}}
\def\bma{\mathbf{a}}
\newcommand{\prtn}{\prt_{\bmn}}
\def\txin{\txt{in}}
\def\txon{\txt{on}}
\maketitle
\section{\bf Introduction}
In this article we consider  equations of the type
\begin{equation}\label{1}
\BAL
-\Gd u+g\circ u&=\mu \ && \text {in} \q \Gw,\\
 u&=\nu \ && \text {on} \q\prt\Gw,
\EAL
\end{equation}
where $\Gw$ is a bounded $C^2$ domain in $\BBR^N$, $g$ is defined
on $\Gw\times\BBR$ and $g\circ u(x)=g(x,u(x))$.

We assume that the nonlinearity $g$ satisfies the conditions,
\begin{equation}\label{2}\BAL
&(a)\q g(x,\cdot)\in C(\BBR),\ g(x,0)=0, \\
&(b)\q g(x,\cdot) \ \text{is a non decreasing and odd function}, \\
&(c)\q g(\cdot,t)\in L^1(\Gw,\rho)  \ \ \forall \ t\in\BBR
\EAL
\end{equation}
where $$\gr(x)=\dist(x,\bdw).$$
The family of functions satisfying these  conditions will be denoted
by $\mathcal{G}_0=\mathcal{G}_0(\Gw)$.

With respect to the data, we assume that $\nu\in\CM(\prt\Gw)$ and $\mu\in\CM(\Gw,\rho)$ where:\vskip 2mm
\noindent $\CM(\bdw)$ denotes the space of bounded Borel measure on $\bdw$ with the usual total variation norm.\vskip 2mm
\noindent $\CM(\Gw,\rho)$ denotes the space of signed Radon measure
$\mu$ in $\Gw$ such that
$$\rho\mu\in\CM(\Gw)\q\text{where}\q
\rho(x):=\text{dist}\,(x,\prt\Gw).$$
\noindent The norm of a measure $\mu\in\CM(\Gw,\rho)$ is given by
$$||\mu||_{\Gw,\rho}=\int_{\Gw}\rho \ d|\mu|.$$
\noindent $L^1(\Gw,\rho)$ denotes the weighted Lebesgue space with weight $\rho$.

We say that $u\in L^1(\Gw)$ is a weak solution of \eqref{1} if $g\circ u\in
L^1(\Gw,\rho)$
and u satisfies the the following,
\begin{equation}\label{3}
\int_{\Gw}\displaystyle\left(-u\Gd\phi+(g\circ u)\phi\right)dx=\int_{\Gw}\phi d\mu-\int_{\bdw}\frac{\prt\phi}
{\prt {\bf n}}d\nu  \  \  \forall\phi\in C^{2}_{0}(\bar\Gw).
\end{equation}
where
$$C_{0}^{2}(\bar\Gw):=\{\phi\in C^2(\bar\Gw):\phi=0 \q \text{on} \q \prt\Gw\}.$$

\bdef{good-meas} Given $g\in\CG_0$, we denote by ${\CM}^g(\Gw)$ the set of all measures $\mu\in\CM(\Gw,\rho)$ such that the boundary value problem
\begin{equation}\label{eq_int}
-\Gd u+g\circ u=\mu \quad \text{in} \ \Gw; \quad
 u=0 \quad \text{on} \ \bdw
\end{equation}
possesses a  weak solution. If $\mu\in\CM^{g}(\Gw)$, we say that $\mu$ is a good measure in $\Gw$.\\
We denote by $\CM^g(\bdw)$ the set of all measures $\nu\in\CM(\bdw)$ such that the boundary value problem
 \begin{equation}\label{eq_bd}
-\Gd u+g\circ u=0 \quad \text{in} \ \Gw; \quad
 u=\nu \quad \text{on} \ \bdw
\end{equation}
possesses a weak solution. If $\nu\in\CM^g(\bdw)$, we say that $\nu$ is a good measure on $\bdw$.\\
Finally, the set of pairs of measures $(\mu,\nu)\in\CM(\Gw,\rho)\times\CM(\bdw)$ such that \eqref{1} possesses a solution will be denoted by $\CM^{g}(\bar\Gw)$.
\es

In a recent paper Marcus and Ponce \cite{MP} studied the following problem. Let $\{\mu_n\}\sbs \CM^g(\Gw)$ be a weakly convergent \seq:  $\mu_n\rightharpoonup\mu$ relative to $C_0(\bar\Gw)$. Let $u_n$ be the solution of \req{eq_int} with $\mu=\mu_n$ and assume that $u_n\to u$ in $L^1(\Gw)$. In general $u$ does not satisfy \req{eq_int}. But it was shown  that there exists a measure $\mu^\#$ \sth
\begin{equation}\label{mu^1}
  -\Gd u+g\circ u=\mu^\#.
\end{equation}
Moreover this measure depends only on the fact that $u_n$ satisfies the equation

 $$-\Gd u_n+g\circ u_n=\mu_n.$$
 It is independent of the boundary data for $u_n$ or indeed on whether $u_n$ has a measure boundary trace. If $\mu_n\geq 0$ then $$0\leq \mu^\#\leq \mu.$$
Furthermore it was shown that, under a mild additional condition on $g$, the following result holds:

\emph{Let $v_n$ be the solution of the Dirichlet problem for equation $-\Gd v=\mu_n$. Suppose that $\mu_n\geq 0$ and that $\{g\circ v_n\}$ is bounded in $L^1(\Gw;\gr)$. Then $\mu^\#$ and $\mu$ are mutually a.c.}

The measure $\mu^\#$ was called the \emph{reduced limit of $\{\mu_n\}$.} This notion is in some sense related to the notion of `reduced measure' introduced in \cite{BMP}. For a specific choice of $\{\mu_n\}$ the reduced limit $\mu^\#$ coincides with the reduced measure. However in general they are not equal.

In the present paper we continue this study considering similar questions \wrto sequences of  pairs $\{(\mu_n,\nu_n)\}\sbs \CM^g(\bar\Gw)$. Such a \seq is called a $g$ good sequence.
 Suppose that $\nu_n\rightharpoonup\nu$ (weak convergence in $\CM(\bdw)$) and that
 $\gr\mu_n\rightharpoonup\tau$ (weak convergence in $\CM(\bar\Gw)$). We shall say that $(\tau,\nu)$ is the weak limit of the \seq.

Let $u_n$ be the solution of \req{1} with $(\mu,\nu)=(\mu_n,\nu_n)$ and suppose that $u_n\to u$ in $L^1(\Gw)$. By \cite{MP}, $u$ satisfies equation \req{mu^1}. Here we show that there exists a measure $\nu^*\in \CM(\bdw)$ \sth $u$ is the weak solution of the \bvp
$$\BAL -\Gd u+g\circ u&=\mu^\# \q \text{in }\Gw\\ u&=\nu^*\q\text{on $\bdw$.}\EAL$$
The pair $(\mu^\#,\nu^*)$ is called the \emph{reduced limit of $\{(\mu_n,\nu_n)\}$.}

In general $\nu^*$  depends on the sequence of pairs  $\{(\mu_n,\nu_n)\}$, not only on $\{\nu_n\}$. If $g$ is subcritical we show that
$$\nu^*=\nu+\tau\chr{\bdw}.$$
However in general the dependence of $\nu^*$ on the sequence of pairs is much more complex.

Here  are some of our main results.

\bth{intr1} Suppose that $\mu_n\geq 0$ and $\nu_n\geq 0$ for every $n\geq1$. If $\nu^\#$ is the reduced limit of the \seq $\{(0,\nu_n)\}$ then
$$0\leq \nu^\#\leq \nu \q\text{and}\q 0\leq\nu^*\leq \nu^\#+\tau\chr{\bdw} .$$
\es

\bth{intr2} Let $v_n$ denote the solution of
$$-\Gd v=\mu_n\q\text{in }\Gw,\qq v=\nu_n\q\text{on }\bdw.$$
Suppose that:

(i) $\mu_n\geq 0$ and $\nu_n\geq 0$ for every $n\geq1$,

 (ii) $\{g\circ v_n\}$ is bounded in $L^1(\Gw;\gr)$,

 (iii) $g\in \CG_0$ satisfies the condition
 $$  \lim_{a,t\tin}\frac{g(x,at)}{ag(x,t)}=\infty\q \text{uniformly \wrto }\; x\in \Gw.$$
 Then $\nu+\tau\chr{\bdw}$ and $\nu^*$ are mutually a.c. In particular $\nu$ and $\nu^\#$ are mutually a.c.
 \es

For the statement of the next result we need an additional definition.
\bdef{negligible}  A nonnegative measure $\gs\in\CM(\bdw)$ is $g$-negligible if
$$\{\gl\in\CM(\bdw):\,0<\gl\leq\gs\}\cap\CM^g(\bdw)=\ems.$$
\es

\bth{intr3}
Assume that $g\in\CG_0$ is convex and satisfies the $\Gd_2$ condition.
Let $\{(\mu_n,\nu_n)\}$ and $\{(\tl \mu_n,\nu_n)\}$ be  $g$-good sequences with weak limits $(\tau,\nu)$ and $(\tl\tau,\nu)$ respectively. Assume that, for every $n\geq 1$, $(|\mu_n|,|\nu_n|)$ and $(|\tl\mu_n|,|\nu_n|)$ are in  $\CM^g(\bar\Gw)$.

Let $u_n$ (resp $\tl u_n$) be the solution of \req{1} with $(\mu,\nu)=(\mu_n,\nu_n)$ (resp. $(\tl\mu_n,\nu_n)$). Assume that
$$u_n\to u,\q \tl u_n\to\tl u \q \text{in }L^1(\Gw)$$
and let $(\mu^*,\nu^*)$ and $(\tl\mu^*,\tl\nu^*)$ denote the reduced limits of $\{(\mu_n,\nu_n)\}$ and $\{(\tl\mu_n,\nu_n)\}$ respectively.

If a \sseq of $\{\gr|\tl\mu_n-\mu_n|\}$ converges weakly in $\CM(\bar\Gw)$ to a measure $\Gl$ \sth $\Gl\chr{\bdw}$ is negligible
then
$$\nu^*=\tl\nu^*.$$
\es
\section{Definitions and auxilliary results.}
\bdef{delta2-cond} Let $g\in\CG_0$. We say that $g$ satisfies the $\Gd_2$ condition if there exists a constant $c>0$ such that
$$g(x,a+b)\leq c\big(g(x,a)+g(x,b)\big)\q\forall\ x\in\Gw,\q a>0,\q b>0.$$
\es

In the next proposition we gather some classical results concerning the \bvp \req{1}.
\bprop{class}
Suppose that $g\in\CG_0,\ \nu\in \CM(\bdw)$ and $\mu\in \CM(\Gw;\gr)$. Then:

\noindent$(i)$ If  $\mu,\nu\geq 0$ then $(\mu,\nu)\in\CM^{g}(\bar\Gw)\Longrightarrow \mu\in\CM^g(\Gw)$ and $\nu\in\CM(\bdw)$.\vskip 1mm

\noindent$(ii)$ If $\mu,\nu\geq 0$, $g$ satisfies $\Gd_2$ condition then
 $$\mu\in\CM^g(\Gw)\q\text {and}\q \nu\in\CM(\bdw)\Longrightarrow (\mu,\nu)\in\CM^{g}(\bar\Gw).$$\vskip 1mm

\noindent$(iii)$ Assume that $(\mu,\nu)\in \CM^{g}(\bar\Gw)$. Then \req{1} possesses a unique solution $u$. This solution satisfies:
\begin{equation}\label{1,est1}
 ||u||_{L^1(\Gw)}+||g\circ u||_{L^{1}(\Gw,\rho)}\leq C(\|\mu\|_{\CM(\Gw;\gr)}+||\nu||_{\CM(\prt\Gw)}),
\end{equation}
where $C$ is a constant depending only on $\Gw$.\vskip 1mm

\noindent$(iv)$ Under the assumption of part (ii), $u\in L^p(\Gw)$ for $1\leq p<\frac{N}{N-1}$ and there exists a constant $C(p)$ depending only on $p$ and $\Gw$ \sth
\begin{equation}\label{1,est2}
||u||_{L^p(\Gw)}\leq C(p)\big(\|\mu\|_{\CM(\Gw;\gr)}+||\nu||_{\CM(\prt\Gw)}\big).
\end{equation}

Furthermore, $u\in W^{1,p}_{loc}(\Gw)$ and for every domain $\Gw'\Subset\Gw$, there exists a constant $C(p,\Gw')$ depending on $p$, $\Gw'$ and $\Gw$ \sth
\begin{equation}\label{1,est3}
||u||_{W^{1,p}(\Gw^{'})}\leq C(p,\Gw')\big(\|\mu\|_{\CM(\Gw')}+||\nu||_{\CM(\prt\Gw)}\big).
\end{equation}
\es
Assertion (i) and (ii) are  obvious (see e.g. \cite[Theorem 2.4.5]{MV}).

Assertion (iii) is due to Brezis \cite{Br}; a proof can be found in \cite{MV} or \cite{LVbook}.

Assertion (iv) is a consequence of (ii) and classical estimates for the corresponding linear problem.

\bdef{tight} A sequence $\{\mu_n\}\in\CM(\Gw,\rho)$ is tight if for every $\epsilon>0$, there exists a neighborhood $U$ of $\bdw$ such that
$$\int_{U\cap\Gw}\rho d|\mu_n|\leq \epsilon.$$
\es

\bdef{exhaustion}A sequence $\{\Gw_n\}$ is an exhaustion of $\Gw$ if $\bar\Gw_n\subset\Gw_{n+1}$ and $\Gw_n\uar\Gw$. We say that an exhaustion $\{\Gw_n\}$ is of class $C^2$ if each domain $\Gw_n$ is of this class. If, in addition, $\Gw$ is a $C^2$ domain and the sequence of domains $\{\Gw_n\}$
is uniformly of class $C^2$, we say that $\{\Gw_n\}$ is a uniform $C^2$ exhaustion.
\es

\bdef{trace} Let $u\in W^{1,p}_{loc}(\Gw)$ for some $p>1$. We say that $u$ possesses an $M-$boundary trace on $\bdw$ if there exists $\nu\in\CM(\bdw)$ such that, for every uniform $C^2$ exhaustion $\{\Gw_n\}$ and every $h\in C(\bar\Gw)$,
$$\int_{\bdw_n}u|_{\bdw_n}h \ dS\to\int_{\bdw}h \ d\nu $$
where $u|_{\bdw_n}$ denotes the Sobolev trace, $dS=d\BBH^{N-1}$ and $\BBH^{N-1}$ denotes $(N-1)$ dimensional Hausdroff measure. The M-boundary trace $\nu$ is denoted by $\tr \,\nu$.
\es

\bprop{th-trace} Let $\mu\in\CM(\Gw,\rho)$ and $\nu\in\CM(\bdw)$. Then a function
$u\in L^1(\Gw)$, with $g\circ u\in L^1(\Gw,\rho)$, satisfies \eqref{3} if and only if
\begin{equation*}
\BAL
-\Gd u+g\circ u &=\mu \q \text{(in the sense of distribution) },\\
\text{tr} \ u &=\nu \q \text{(in the sense of Definition \ref{d:trace}) }.
\EAL
\end{equation*}
\es

This is an immediate consequence of
\cite[Proposition 1.3.7]{MV}.

\section{\bf Reduced limit of a \seq of measures in $\CM^g(\bdw)$}
In this section we discuss a sequence of problems
\begin{equation}\label{bvpn}
\BAL
-\Gd u_n+g\circ u_n&=0 && \text {in } \Gw,\\
u_n&=\nu_n && \text {on } \prt\Gw,
\EAL
\end{equation}
where $g\in\CG_0$ and
\begin{equation}\label{nun}\BAL
&(i)\quad   \nu_n\;\text{is a  good measures on } \bdw \forevery n\in \BBN\\
&(ii)\quad\nu_n\rightharpoonup\nu \;\text{in } \CM(\bdw)
\EAL
\end{equation}

\blemma{2.1} Assume that $\{\nu_n\}$ satisfies \req{nun} and let $u_n$ be the solution of \req{bvpn}.
Then there exists a subsequence $\{u_{n_k}\}$ that converges in $L^1(\Gw)$.
\es
\begin{proof} By \rprop{class}, $\{u_n\}$ is bounded in $L^p(\Gw)$ for every $p\in [1,\,N/(N-1))$.
\Consy,  for each such $p$, $\{u_n\}$ is uniformly integrable in $L^p(\Gw)$. Furthermore $\{u_n\}$ is bounded in $W^{1,p}\loc(\Gw)$ for every $p\in [1,\,N/(N-1))$. Therefore there is a subsequence
$\{u_{n_k}\}$ which converges in $L^1\loc(\Gw)$ and pointwise a.e. to a function $u$. Combining these facts we conclude that $u_{n_k}\to u$ in $L^1(\Gw)$.
\end{proof}

\medskip

To simplify the presentation we introduce the following:
 \bdef{bar-con} (i) Let $\{\mu_n\}$ be a bounded \seq of measures in $\CM(\Gw;\gr)$.
Assume that $\gr\mu_n$ is extended to a Borel measure $(\mu_n)_\gr\in \CM(\bar \Gw)$ defined as zero on $\bdw$. We say that $\{\gr\mu_n\}$ \emph{converges weakly in $\bar\Gw$} to a measure $\tau\in \CM(\bar\Gw)$
if $\{(\mu_n)_\gr\}$ converges weakly to $\tau$ in $\CM(\bar\Gw)$, i.e.
$$\int_\Gw \gf \gr\,d\mu_n\to\int_{\bar\Gw} \gf \,d\tau \forevery \gf\in C(\bar\Gw).$$
 This convergence is denoted by
$$\gr\mu_n\underset{\bar\Gw}{\rightharpoonup}\tau.$$
(ii) Let $\{\mu_n\}$ be a \seq in  $\CM_{loc}(\Gw)$. We say that the \seq converges weakly to $\mu\in \CM_{loc}(\Gw)$ if it converges in the distribution sense, i.e.,
$$\int_\Gw \gf\,d\mu_n\to \int_\Gw \gf\,d\mu \forevery \gf\in C_c(\Gw).$$
This convergence is denoted by $\mu_n\underset{d}{\rightharpoonup}\mu$.
\es

 If $\{\gr\mu_n\}$ converges weakly in $\bar\Gw$ to $\tau$ then $$\mu_n\underset{d}{\rightharpoonup}\mu_{int}:=\frac{\tau}{\gr}\chr{\Gw}.$$
Thus, for $\tau $  as in part (i),
\begin{equation}\label{tau_bd}
   \tau=\tau\chr{\bdw}+\gr\mu_{int}.
\end{equation}

\blemma{prtn}
Let $\{\mu_n\}$ be as in \rdef{bar-con}(i) and assume that $\gr\mu_n\underset{\bar\Gw}{\rightharpoonup}\tau.$ Then
\begin{equation}\label{prtn}
 \lim_{n\tin}  \int_\Gw \vgf\,d\mu_n =\int_\Gw\vgf\,d\mu_{int}-\int_{\bdw}\frac{\prt\vgf}{\prt \bmn}d\tau
\end{equation}
for every $\vgf\in C_0^1(\bar\Gw)$.
\es
\begin{proof} Put
\begin{equation}\label{bar-phi}
\bar\vgf=\begin{cases}\vgf /\rho&\text{in }\Gw\\ -\frac{\prt \vgf}{\prt\mathbf{n}}&\text{on }\bdw.\end{cases}
\end{equation}
Then $\bar\vgf\in C(\bar \Gw)$ and \consy, using \req{tau_bd},
$$\BAL
 \lim_{n\tin}\int_\Gw \vgf\,d\mu_n&=\lim_{n\tin}\int_\Gw \bar\vgf \gr\,d\mu_n=\int_{\bar\Gw} \bar\vgf \,d\tau\\
 &= \int_\Gw\vgf\,d\mu_{int}-\int_{\bdw}\frac{\prt\vgf}{\prt \bmn}d\tau.
\EAL$$
\end{proof}

\bth{th1}  Assume that $g\in \CG_0$ and that $\{\nu_n\}$ is a \seq of measures satisfying \req{nun}.
Let $u_n$ be the solution of \req{bvpn} and assume  that
\begin{equation}\label{untou}
 u_n\to u\q\text{in }L^1(\Gw).
\end{equation}
 Then there exists a measure $\nu^\#\in \CM^g(\bdw)$ \sth
\begin{equation}\label{lim_4}
\BAL
-\Gd u+g\circ u&=0 && \text {in } \Gw,\\
u&=\nu^\# && \text {on } \prt\Gw.
\EAL
\end{equation}

Furthermore $\{(g\circ u_n)\rho\}$ converges weakly in $\bar\Gw$ to a measure $\gl\in \CM(\bar\Gw)$ and
 \begin{equation}\label{lim_nun}
\nu^\#=\nu-\gl\chr{\bdw}.
\end{equation}
If $\nu_n\geq 0$ then $0\leq \nu^\#\leq \nu$.
\es

{\Remark}  The measure $\nu^{\#}$ defined above is called \emph{the reduced limit of $\{\nu_n\}$}. We emphasize that $\nu^\#$ depends on the \seq, not only on its limit.

\begin{proof} By assumption $-\Gd u_n +g\circ u_n=0$ in $\Gw$ and $u_n\to u$ in $L^1(\Gw)$ . Therefore, by \cite[Thm. 1.3]{MP},
\begin{equation}\label{t1.1}
   -\Gd u+g\circ u=0 \q\text{in }\Gw.
\end{equation}
(Note that,   in the notation of \cite{MP},  the present case corresponds to $\mu_n=0$ and therefore $\mu^\#=0$.)

Consider a \sseq of $\{u_n\}$ \sth $\{\gr\,g\circ u_n\}$ converges weakly in $\bar \Gw$. The subsequence is still denoted by $\{u_n\}$ and we denote by $\gl$ the weak limit of $\{\gr\,g\circ u_n\}$ in $\CM(\bar\Gw)$. Put
\begin{equation}\label{gl1}
\gl_{in}=\gl\chr{\Gw},\q\gl_{bd}=\gl\chr{\bdw}.
\end{equation}
Then  $g\circ u_n\underset{d}{\rightharpoonup} \gl_{in}/\gr$ and \consy
$$-\Gd u+\frac{\gl_{in}}{\rho}=0 \q\text{in }\Gw.$$
Comparing with \req{t1.1} we obtain,
\begin{equation}\label{glin}
   \gl_{in}=\rho(g\circ u).
\end{equation}



For every $\vgf\in C^2_0(\bar\Gw)$,
\begin{equation}\label{vgf}
-\int_\Gw\, u_n\Gd \vgf \, dx +\int_\Gw (g\circ u_n)\vgf\,dx = -\int_{\bdw} \frac{\prt\vgf}{\prt\mathbf{n}}d\nu_n.
\end{equation}
\bcom
Note that the function
\begin{equation}\label{bar-phi}
\bar\vgf=\begin{cases}\vgf /\rho&\text{in }\Gw\\ -\frac{\prt \vgf}{\prt\mathbf{n}}&\text{on }\bdw\end{cases}
\end{equation}
is in $C^1(\bar \Gw)$.
\end{comment}
By the definition of $\gl$, \rlemma{prtn} and \req{glin},
\begin{equation}\label{weak-phi}
\lim_{n\tin}\int_\Gw (g\circ u_n)\vgf\,dx =\int_\Gw (g\circ u)\vgf\,dx-\int_{\bdw}\frac{\prt\vgf}{\prt\bmn}d\gl
\end{equation}
Therefore, taking the limit in \req{vgf} we obtain
$$-\int_\Gw\, u\Gd \vgf \, dx +\int_\Gw (g\circ u)\vgf\,dx -\int_{\bdw}\frac{\prt\vgf}{\prt\bmn}d\gl=- \int_{\bdw} \frac{\prt\vgf}{\prt\mathbf{n}}d\nu.$$
Thus $u$ is a weak solution of \req{lim_4} with $\nu^\#$ given by \req{lim_nun}.
By \rprop{th-trace}, $\nu^\#$ is the  M-boundary trace of $u$; hence $\nu^\#$ is independent of the specific subsequence of $\{(g\circ u_n)\rho\}$ that converges weakly in $\bar\Gw$. This fact and \req{lim_nun} imply that $\gl_{bd}$ is independent of the \sseq. By \req{glin}, $\gl_{in}$ is independent of the \sseq. Therefore the full sequence $\{\rho(g\circ u_n)\}$ converges to $\gl$.

If $\nu_n\geq 0$ then $u_n\geq 0$ and $g\circ u_n\geq 0$. Therefore, in this case, $\gl\geq0$ and \consy $\nu^\#\leq \nu$. Further, $u\geq0$ and therefore its M-boundary trace, namely  $\nu^\#$, is non-negative.
\end{proof}

\blemma{comp1} Let $\{\nu_n\}$ and $\{\nu'_n\}$ be sequences of measures in $\CM^g(\bdw)$ with weak limits $\nu$ and $\nu'$ respectively. Let $u_n$ (resp. $u'_n$) be the solution of \req{1} with $\mu=0$ and $\nu=\nu_n$ (resp. $\nu=\nu'_n$). Assume that $u_n\to u$ and $u'_n\to u'$ in $L^1(\Gw)$.

If $\nu_n\leq \nu'_n$ for every $n$ then $\nu^\#$ and $(\nu')^\#$ (the reduced limits of the two sequences)
satisfy
\begin{equation}\label{comp1}
 0\leq (\nu')^\#-\nu^\#\leq \nu'-\nu.
\end{equation}
\es
\bproof Since $\nu_n\leq \nu'_n$ we have $u_n\leq u'_n$. Hence
$$\nu^\#=\tr\,u\leq \tr\,u'=(\nu')^\#$$
and
$$\gl=\lim\rho (g\circ u_n)\leq \lim\rho (g\circ u'_n)=\gl'.$$
By \rth{th1} these limits exist in the sense of weak convergence in $\CM(\bar\Gw)$. Furthermore,
$$\nu^\#=\nu-\gl_{bd},\q (\nu')^\#=\nu'-\gl'_{bd}.$$
Hence
$$(\nu')^\#-\nu^\#=(\nu'-\nu)-(\gl'_{bd}-\gl_{bd})\leq \nu'-\nu.$$
\eproof

\bth{th2} In addition to the assumptions of \rth{th1}, assume that $g$ satisfies
\begin{equation}\label{g_at_infty}
   \lim_{a,t\tin}\frac{g(x,at)}{ag(x,t)}=\infty\q \text{uniformly \wrto }\; x\in \Gw.
\end{equation}
Put $v_n:=\BBP(\nu_n)$, i.e.
\begin{equation}\label{Pnun}
  -\Gd v_n=0 \q \text{in }\Gw,\q v_n=\nu_n\q\text{on }\bdw.
\end{equation}

If $\nu_n\geq0$ and $\{g\circ v_n\}$ is bounded in $L^1(\Gw;\rho)$ then $\nu$ and $\nu^\#$ (the reduced limit of $\{\nu_n\}$) are mutually absolutely continuous.
\es

We postpone the proof to Section 3 where we present a more general version of this result.

\bprop{nu_pm} Assume that $g\in \CG_0$.
Let $\{\nu_n\}\subset\CM(\bdw)$ be a bounded sequence \sth $|\nu_n|\in \CM^g(\bdw)$ for every $n$.
Denote by $u_n$, $u_{n,1}$ and $u_{n,2}$ the solution of \req{1} with $\mu=0$ and $\nu=\nu_n$, $\nu=\nu_n^+$
and $\nu=-\nu_n^-$ respectively.
Assume that
\begin{equation}\label{nu+,-}
\nu_n^{+}\warrow \nu^{+} \ \ \text{and}\ \ \nu_n^{-}\warrow \nu^{-} \ \ \text{in} \ \ \CM(\bdw).
\end{equation}

Then $\{u_n\}$ converges in $L^1(\Gw)$ if and only if  $\{u_{n,1}\}$ and $\{u_{n,2}\}$ converge in $L^1(\Gw)$. Assuming the convergence of these sequences,
denote by $\nu^\#$, $\nu_1^\#$ and $\nu_2^\#$ the reduced limits of  $\{\nu_n\}$, $\{\nu_n^+\}$ and $\{-\nu_n^-\}$ respectively.
Then
\begin{equation}\label{nu1,2}
  \nu_{1}^{\#}=(\nu^{\#})^{+}, \q\nu_{2}^{\#}=-(\nu^{\#})^{-}.
\end{equation}

In particular $$\nu^{\#}=\nu_{1}^{\#}+\nu_{2}^{\#}$$ and
$$\nu^{\#}=\nu \ \ \text{if and only if} \ \ \nu_{1}^{\#}=\nu^{+} \ \ \text{and} \ \ \nu_{2}^{\#}=-\nu^{-}.$$
\es
\begin{proof} First assume that $\{u_n\}$, $\{u_{n,1}\}$ and $\{u_{n,2}\}$ converge in $L^1(\Gw)$. In that case, \req{nu1,2} is proved exactly in the same way as  \cite[Proposition 7.3]{MP},
 using Lemma \ref{l:comp1} and the last assertion of  \rth{th1}.

Next assume that $\{u_n\}$ converges in $L^1(\Gw)$ and let $\nu^\#$ be the reduced limit of $\{\nu_n\}$.  Extract a \sseq $\{u_{n_k}\}$ \sth  $\{u_{n_k}^+\}$ and $\{u_{n_k}^-\}$ converge in $L^1(\Gw)$. Denote the limits of these \seq{s} by $u'$ and $u''$ respectively. By \req{nu1,2}
$$\tr\,u'=\nu_1^\#=(\nu^\#)^+.$$
Thus $u'$ is independent of the \sseq previously extracted. This implies that $u_n^+\to u'$ in $L^1(\Gw)$.
Similarly we conclude that $u_n^-\to u''$ in $L^1(\Gw)$.

The same argument shows that if $\{u_{n,1}\}$ and $\{u_{n,2}\}$ converge in $L^1(\Gw)$ then $\{u_n\}$ converges in $L^1(\Gw)$.
\end{proof}

As a consequence of this proposition one obtains the following extension of \rth{th2} to sequences of signed measures.

\bcor{th2'} In addition to the assumptions of \rprop{nu_pm} assume that $g$ satisfies \req{g_at_infty}. Let
$\bar v_n=\BBP(|\nu_n|)$ and assume that $\{g\circ\bar v_n\}$ is bounded in $L^1(\Gw;\gr)$. Then $\nu^\#$ and $\nu$ are mutually absolutely continuous. More precisely, $(\nu^\#)^+$ and $\nu^+ $ (respectively  $(\nu^\#)^-$ and $\nu^- $) are mutually a.c.
\es

\section{\bf Reduced limit of a sequence of pairs in $\CM^g(\bar \Gw)$}
In this section we discuss the reduced limit of a sequence of pairs $\{(\mu_n,\nu_n)\}\subset\CM^g(\bar\Gw)$ associated with problem,
\begin{equation}\label{bvpn1}
\BAL
\Gd u_n+g\circ u_n &=\mu_n \q &&\text{in}\q \Gw, \\
u_n&=\nu_n \q &&\text{on}\q \bdw.
\EAL
\end{equation}
We assume that $\nu_n$ satisfies \eqref{nun} and $\mu_n$ satisfies
\begin{equation}\label{nun2}\BAL
&(i)\q\mu_n\ \text{is a good measure in}\ \Gw\q\forall\ n\in\BBN \\
&(ii)\q\rho\mu_n\underset{\bar\Gw}{\rightharpoonup}\tau\in \CM(\bar\Gw)\q \text{(see \rdef{bar-con})}.
\EAL
\end{equation}

\bth{th4}
Assume that $g\in\CG_0$, $(\mu_n,\nu_n)\in \CM^g(\bar\Gw)$, $\{\nu_n\}$ satisfies \eqref{nun} and $\{\mu_n\}$ satisfies \eqref{nun2}. Let $u_n$ be the solution of \eqref{bvpn1} and assume that
$$u_n\to u\q\text{in} \ L^1(\Gw).$$
Then:

\nind $(i)\q \{\gr(g\circ u_n)\}$ converges weakly in $\bar \Gw$ and

\nind$(ii)\q \exists\ \mu^{*}\in\CM(\Gw,\rho), \;\nu^{*}\in\CM(\bdw)$ such that
\begin{equation}\label{mn-star}\BAL
-\Gd u+g\circ u &=\mu^{*} \q &&\text{in}\q\Gw\\
u&=\nu^{*}\q&&\text{on}\q\bdw.
\EAL
\end{equation}

Furthermore, if $\mu_n\geq0$ and $\nu_n\geq 0$ for every $n\geq1$ then
\begin{equation}\label{ineq_nu,1}
   0\leq\nu^{*}\leq(\nu+\tau\chr{\bdw}).
\end{equation}

\es
{\bf \Remark} By \cite[Theorem 1.3]{MP}, $\mu^{*}$ is independent of $\nu_n$.

\begin{proof}
Our assumptions imply that  $\{\nu_n\}$ is bounded in $\CM(\bdw)$ and $\{\mu_n\}$ is bounded in $\CM(\Gw;\gr)$. Hence  $\{\rho(g\circ u_n)\}$ is bounded in $L^1(\Gw)$. Therefore there exists a subsequence (still denoted by $\{u_n\}$) \sth
$$\gr\,g\circ u_n\underset{\bar\Gw}{\rightharpoonup}\gl$$
(see \rdef{bar-con}).
Put
$$ \gl_{int}=\frac{\gl}{\gr}\chr{\Gw} \ \text{and}\ \gl_{bd}=\gl\chr{bd}.$$
By \rlemma{prtn},
\begin{equation}\label{glint}
   \lim_{n\tin}\int_\Gw (g\circ u_n)\vgf\,dx=\int_\Gw \vgf\,d\gl_{int}- \int_{\bdw}\frac{\prt\vgf}{\prt\bmn}d\gl
\end{equation}
and \req{prtn} holds for every $\vgf\in C^2_0(\bar\Gw)$.

As $u_n$ is the weak solution of \req{bvpn1},
\begin{equation}\label{22}
\int_{\Gw}\big(-u_n\Gd\vgf+(g\circ u_n)\vgf\big)dx=\int_{\Gw}\vgf\ d\mu_n-\int_{\bdw}\frac{\prt\vgf}{\prt \bf n} d\nu_n
\end{equation}
for every $\vgf\in C^2_0(\bar\Gw)$.
Taking the limit as $n\tin$ and using \req{prtn} and \req{glint} we obtain,
$$-\int_{\Gw}u\Gd\vgf dx+\int_{\Gw}\vgf\,d(\gl_{int}-\mu_{int})-\int_{\bdw}\frac{\prt\vgf}{\prt\bf n}d(\gl_{bd}-\tau_{bd})=-\int_{\bdw}\frac{\prt\vgf}{\prt\bf n}d\nu$$
for every $\vgf\in C^2_0(\bar\Gw)$. Thus $u$ is the  weak solution of \req{mn-star} where
\begin{equation}\label{3.1}
  \mu^*=g\circ u-(\gl_{int}-\mu_{int})
\end{equation}
\begin{equation}\label{3.2}
    \nu^*:=\nu-(\gl_{bd}-\tau_{bd}).
\end{equation}

By \cite[Theorem 1.3]{MP},
$\mu^*$ depends on $\{\mu_n\}$ but is independent of $\{\nu_n\}$.

The fact that $u$  is the weak solution of \req{mn-star} implies that $\nu^{*}$ is the M-boundary trace of $u$; as such $\nu^{*}$ is independent of the specific weakly convergent subsequence of $\{\rho(g\circ u_n)\}$. Therefore, by \req{3.2}, $\gl_{bd}$ is independent of the subsequence. In addition by \req{3.1} and \cite[Theorem 1.3]{MP},   $\gl_{int}$ is independent of the subsequence. This implies that the full sequence $\{\rho(g\circ u_n)\}$ converges to $\gl$.

If $\mu_n,\ \nu_n\geq 0$ then $u_n\geq 0$ and $g\circ u_n\geq 0$. Therefore, in this case, $\nu^*\geq0$ and $\gl\geq0$; hence, by \req{3.2} $\nu^{*}\leq \nu+\tau_{bd}$.
\end{proof}

\bdef{mun,nun} If $\{(\mu_n,\nu_n)\}\in\CM^g(\bar\Gw)$, $\{\nu_n\}$ satisfies \eqref{nun} and $\{\mu_n\}$ satisfies \eqref{nun2} we say that $\{(\mu_n,\nu_n)\}$ is a $g$-good \seq  that converges weakly to  $(\tau,\nu)$ in $\bar\Gw$.

If in addition $u_n\to u$ in $L^1(\Gw)$ we say that $(\mu^*,\nu^*)$, defined as in \rth{th4}, is the reduced limit and $\nu^*$ is the boundary reduced limit of $\{(\mu_n,\nu_n)\}$.
\es

\bth{th6}In addition to the assumptions of Theorem \ref{t:th4}, assume that $g(x,\cdot)$ satisfies \req{g_at_infty}.

Let $v_n$ be weak solution of
\begin{equation}\label{Pnun'}
  -\Gd v_n=\mu_n \q \text{in }\Gw,\q v_n=\nu_n\q\text{on }\bdw.
\end{equation}
\nind If $\mu_n,\ \nu_n\geq0$ and $\{g\circ v_n\}$ is bounded in $L^1(\Gw;\rho)$ then
$\nu^{*}$ (defined as in Theorem \ref{t:th4}) and $\nu+\mu_{bd}$ are absolutely continuous with respect to each other.
\es
\Remark As a consequence of  \cite[Theorem 8.1]{MP} in combination with \cite[Theorem 1.3]{MP}, $\mu^{*}$ and $\mu_{int}$ are absolutely continuous \wrto each other.

\begin{proof} Given $\ga\in(0,1)$, we have $$0\leq g\circ(\ga v_n)\leq g\circ v_n.$$
Thus there exists $C_0>0$ such that
$$\|g\circ(\ga v_n)\|_{L^1(\Gw,\rho)}\leq C_0   \forevery n\geq 1, \;\forall\ga\in (0,1).$$
Let $\{\ga_k\}$ be a \seq decreasing to zero. One can extract a \sseq of $\{\rho(g\circ(\ga v_n))\}$ (still denoted $\{\rho(g\circ(\ga v_n))\}$) \sth, for each $k$, there exists a measure $\gs_k\in \CM(\bar\Gw)$ \sth
\begin{equation}\label{gsk}
  \gr g\circ (\ga_k\, v_n)\underset{\bar\Gw}{\rightharpoonup}\gs_k.
\end{equation}

\bcom
First we prove,
\begin{equation}\label{claim1}
 \ga_k(\nu+\tau\chr{\bdw})\leq \gs_k\chr{\bdw}+\nu^{*}.
\end{equation}
\end{comment}

Let $w_{n,k}$ be the solution of the problem
\begin{equation}\label{13}
\BAL
-\Gd\,w+g\circ w &= \ga_k\mu_n\q \text {in} \;\Gw,\\
  w&=\ga_k\nu_n \q \text {on}\;\bdw.
\EAL
\end{equation}

$\ga_k v_n$ is a supersolution of problem \req{13}; therefore
\begin{equation}\label{gav}
w_{n,k}\leq   \ga_k\,v_n.
\end{equation}

Passing to a subsequence, we may assume that $\{w_{n,k}\}$ converges in $L^1(\Gw)$ for each $k\in \BBN$. Denote by
$(\mu_k^{*},\nu_{k}^{*})$ the reduced limit of $(\ga_k\mu_n,\ga_k\nu_n)$. By \rth{th4} $\{\gr(g\circ w_{n,k})\}$ converges weakly in $\bar \Gw$ for each $k\in \BBN$; we denote its limit by $\gl_k$. By the proof of \rth{th4}-- specifically \req{3.2} --
$$\nu_k^*=\ga_k\nu -(\gl_k-\ga_k\tau)\chr{\bdw}.$$

By \req{gav}
$$\rho(g\circ(\ga_k v_n)-\rho((g\circ w_{n,k}))\underset{\bar\Gw}{\rightharpoonup}\gs_k-\gl_k\geq0.$$
Thus
\begin{equation}\label{3.0*}
   (\gs_k-\gl_k)\chr{\bdw}=\gs_k\chr{\bdw} +\nu_k^*-\ga_k(\nu+\tau\chr{\bdw})\geq 0.
\end{equation}

Let $u_n$ be the solution of \req{bvpn1}. Evidently $w_{n,k}\leq u_n$ for every $k,n\in \BBN$. \Consy
$$w_k:=\lim w_{n,k}\leq \lim u_n=u.$$
This in turn implies that
\begin{equation}\label{nuk<nu}
   \nu_k^*=\tr\,w_k\leq\tr\,u\leq \nu^*.
\end{equation}
Finally, combining  \req{3.0*} and \req{nuk<nu} we obtain
\begin{equation}\label{3.1*}
    \ga_k(\nu+\tau\chr{\bdw})\leq\gs_k\chr{\bdw}+\nu^{*}.
\end{equation}

In view of \req{g_at_infty}, for every $\epsilon>0$ there exist $a_0, \, t_0>1$, such that
\begin{equation}\label{14}
 \frac{g(x, at)}{ag(x, t)}\geq\frac{1}{\epsilon} \   \  \forall a\geq a_0, \ \ t\geq t_0.
\end{equation}
Consider the splitting of $\rho(g\circ(\ga_k v_n))$ as follows,
\begin{equation*}
\rho(g\circ(\ga_k v_n))=\rho( g\circ(\ga_k v_n))\chr{[\ga_k v_n<t_0]}+
\rho(g\circ(\ga_k v_n))\chr{[\ga_k v_n\geq t_0]}.
\end{equation*}
By passing to a \sseq we may assume that each of the terms on the right hand side converges weakly in $\bar\Gw$ to $\gs_{1,k}$ and $\gs_{2,k}$ respectively, for each $k\geq 1$.
Since $\{\rho( g\circ(\ga_k v_n))\chr{[\ga_k v_n<t_0]}\}$ is uniformly bounded, $\gs_{1,k}\chr{\bdw}=0$. Thus
$$\gs_k\chr{\bdw}=\gs_{2,k}\chr{\bdw}.$$
But
$$\norm{\gs_{2,k}}_{\CM(\bar\Gw)}\leq \liminf_{n\to\infty}\int_{[\ga_k v_n\geq t_0]}\rho( g\circ(\ga_k v_n)).$$
Therefore
\begin{equation*}
||\gs_k\chr{\bdw}||_{\CM(\bdw)}\leq
\liminf_{n\to\infty}\int_{[\ga v_n\geq t_0]}\rho( g\circ(\ga_k v_n)).
\end{equation*}
For $k$ sufficiently large, say $k\geq k_\ge$, $\rec{\ga_k}\geq a_0$ . Applying \eqref{14} with $a=\frac{1}{\ga_k},\ \ t=\ga_k v_n$, we obtain
$$g\circ(\ga_k v_n)\chr{[\ga v_n\geq t_0]}\leq\ga_k\epsilon(g\circ v_n)$$
for $k\geq k_\ge$ and $n\geq 1$.
Hence
\begin{equation*}
||\gs_k\chr{\bdw}||_{\CM(\bdw)} \leq \ga_k\epsilon\liminf_{n\to\infty}\int_{\Gw}\rho(g\circ v_n)\leq C_0\epsilon\ga_k \forevery k\geq k_\ge.
\end{equation*}
Therefore
\begin{equation}\label{3.2*}
  \frac{\norm{\gs_k\chr{\bdw}}_{\CM(\bdw)}}{\ga_k}\to 0.
\end{equation}
\vspace{2mm}

 Since
$\nu^{*}\leq\nu+\tau\chr{\bdw}$, we only have to
prove that $\nu+\tau\chr{\bdw}$ is absolutely continuous \wrto $\nu^{*}$. Let $E\subset\bdw$ be a Borel set \sth
$\nu^{*}(E)=0$. Then, by \req{3.1*}
$$\ga_k(\nu(E)+\tau(E))\leq\gs_k(E) \ \ \forall\ k\geq 1.$$
This inequality and \req{3.2*} imply that $\nu(E)+\tau(E)=0$.
\end{proof}

\blemma {comp2} Let $g\in\CG_0$.
Assume that $\{(\mu_n,\nu_n)\}$ and $\{(\tl \mu_n\},\tl\nu_n)\}$ be $g$  good sequences converging weakly in $\bar\Gw$  to $(\tau,\nu)$ and $(\tl\tau,\tl\nu)$ respectively.

 Let $u_n$ (resp $\tl u_n$) be the solution of \req{1} with $(\mu,\nu)=(\mu_n,\nu_n)$ (resp. $(\tl\mu_n,\tl\nu_n)$). Assume that
$$u_n\to u,\q \tl u_n\to\tl u \q \text{in }L^1(\Gw)$$
and let $(\mu^*,\nu^*)$ and $(\tl\mu^*,\tl\nu^*)$ denote the reduced limits of $\{(\mu_n,\nu_n)\}$ and $\{(\tl\mu_n.\tl\nu_n)\}$ respectively.

Under these assumptions, if
$$\mu_n\leq\tilde\mu_n,\q \nu_n\leq\tilde\nu_n \forevery  n\geq 1$$
 then
 \begin{equation}\label{comp2}\BAL
 &(a)\q 0\leq\tilde\nu^{*}-\nu^{*}\leq(\tl\nu-\nu) +(\tl\tau-\tau)\chr{\bdw},\\
   &(b)\q 0\leq\tilde\mu^{*}-\mu^{*}\leq  \rec{\gr}(\tl\tau-\tau)\chr{\Gw}=:\tl\mu_{int}-\mu_{int} .
 \EAL\end{equation}
 \es
\begin{proof} Inequality \req{comp2} (b) is proved in \cite[Theorem 7.1]{MP}. (Recall that the reduced limit $\mu^*$ is independent of $\{\nu_n\}$.) It remains to prove \req{comp2}(a).
Clearly $u_n\leq\tl u_n$, thus $u\leq\tl u$. Hence $\nu^*\leq\tl\nu^*$.
By \rth{th4}  there exist measures $\gl, \tl\gl\in\CM(\bar\Gw)$ such that
$$\rho g\circ u_n\underset{\bar\Gw}{\warrow}\gl \q\text{and}\q
\rho g\circ\tl u_n\underset{\bar\Gw}{\warrow}\tl\gl.$$ Since $u_n\leq\tl u_n$, we also have $\gl\leq\tl\gl$. Therefore from Theorem \ref{t:th4}
\begin{equation*}
\BAL
0\leq\tl\nu^*-\nu^* &=\tl\nu+\tl\tau\chr{\bdw}-\tl\gl\chr{\bdw}-(\nu+\tau\chr{\bdw}-\gl\chr{\bdw})\\
&=(\tl\nu-\nu)+(\tl\tau-\tau)\chr{\bdw}-(\tl\gl-\gl)\chr{\bdw}\\
&\leq (\tl\nu-\nu)+(\tl\tau-\tau)\chr{\bdw}
\EAL
\end{equation*}
This proves \req{comp2}(a).
\end{proof}

\bcor{tight} Let $g\in\CG_0$, $u_n$ be the weak solution of \eqref{bvpn1} and $v_n$ be the weak solution of \eqref{1} with $(\mu,\nu)=(\tilde\mu_n,\nu_n)$. Assume that
\begin{equation}\label{t_1}
\BAL
\rho\mu_n\underset{\bar\Gw}{\warrow}\mu, \q \rho\tl\mu_n\underset{\bar\Gw}{\warrow}\tl\mu\q\text{and}\q    \nu_n\warrow\nu;\\
u_n\to u, \q v_n\to v \q\text{in}\q L^1(\Gw).
\EAL
\end{equation}
Let $(\mu^*,\nu^*)$  $($respectively $ (\tl\mu^*,\tl\nu^*))$ denote the reduced limit of $\{(\mu_n,\nu_n)\}$ (respectively $\{ (\tl\mu_n,\nu_n)\}$).
If $\mu_n\leq\tl\mu_n$ and $\{\tl\mu_n-\mu_n\}$ is tight then $\nu^*=\tl\nu^{*}$.
\es
\begin{proof}
By Lemma \ref{l:comp2},
$$0\leq\tl\nu^*-\nu^*\leq(\tl\tau-\tau)\chr{\bdw}.$$
Since $\{\tl\mu_n-\mu_n\}$ is tight we have $(\tl\tau-\tau)\chr{\bdw}=0$ and \consy $\nu^*=\tl\nu^{*}$.
\end{proof}

The next corollary provides an improved inequality for $\nu^*$ (compare to \req{ineq_nu,1}).
\bcor{ineq_nu} Let $\{(\mu_n,\nu_n)\}$ be a $g$-good sequence weakly convergent to $(\tau,\nu)$ in $\bar\Gw$
(in the sense of \rdef{mun,nun}). Assume that the sequence has reduced limit $(\mu^*,\nu^*)$.

If $\mu_n\geq0$ and $\nu_n\geq0$ for every $n\geq1$ then
\begin{equation}\label{ineq_nu}
  \nu^\#\leq \nu^*\leq \nu^\#+\tau\chr{\bdw},
\end{equation}
where $\nu^\#$ is the reduced limit of $\{\nu_n\}$ defined in Section 2.
\es
\begin{proof} We apply \rlemma{comp2} to the sequences $\{(\mu_n,\nu_n)\}$ and $\{(0,\nu_n)\}$
\end{proof}

\section{\bf Subcritical problem}


\bth{th5} Assume that $g\in \CG_0$ has subcritical growth with respect to the boundary, i.e., there exists $C>0$ and $q<\frac{N+1}{N-1}$ such that
\begin{equation}\label{subcr0}
   |g(x,t)|\leq C(|t|^q+1)\q \forall t\in\BBR.
\end{equation}

 Let $\{\mu_n\}\sbs\GTM(\Gw;\gr)$ and $\{\nu_n\}\sbs \GTM(\bdw)$ and let $u_n$ be the weak solution of the problem
\begin{equation}\label{26}
\BAL
-\Gd u_n+g\circ u_n &=\mu_n \q &&\text {in} \q \Gw,\\
u_n &=\nu_n \q &&\text {on} \q \bdw.
\EAL
\end{equation}
Assume that
\begin{equation}\label{subcr1}
   \nu_n\rightharpoonup\nu \q\text{weakly in $\bdw$},\q \gr\mu_n\underset{\bar\Gw}{\rightharpoonup}\tau\q\text{weakly in $\bar\Gw$.}
\end{equation}

If
$u_n\to u$ in $L^1(\Gw)$ then $u$ is a weak solution of the problem
\begin{equation}\label{27}
\BAL
-\Gd u+g\circ u &=\mu_{int} \q &&\text {in}\q \Gw,\\
 u &=\nu+\tau\chr{\bdw} \q  &&\text {on}\q \bdw.
\EAL
\end{equation}
where $\mu_{int}=\frac{\tau}{\rho}\chr{\Gw}$.
\es
\Remark In the present case, if $\mu_n$, $\nu_n$ satisfy the assumptions of the theorem then $\{u_n\}$ has a \sseq converging in $L^1(\Gw)$. This is proved as in Section 2.

\nind{\it Notation:} Given $\mu\in \GTM(\Gw;\gr)$ we denote by $\BBG(\mu)$, the weak solution of the problem
\begin{equation}\label{28}
-\Gd u =\mu \q  \text {in} \q \Gw; \q u=0 \q  \text {on} \q \bdw.
\end{equation}
Given $\nu\in \GTM(\bdw)$ we denote by $\BBP(\nu)$ the weak solution of the problem
\begin{equation}\label{28'}
\Gd v =0 \q  \text {in} \q \Gw; \q v=\nu \q  \text {on} \q \bdw.
\end{equation}

\begin{proof} First we show that
\begin{equation}\label{subcr2}
 g\circ u_n\to g\circ u\q\text{in }L^1(\Gw,\rho).
\end{equation}

Define $\mathbb{G}(|\mu_n|):=v_n$ and $\BBP(|\nu_n|):=v'_n$. Then $v_n+v'_n$ satisfies
$$-\Gd(v_n+v'_n)=|\mu_n|\q \text{in}\ \Gw;\q v_n+v'_n=|\nu_n|\q\text{on}\ \bdw.$$ Let $U_n$ denote the weak solution of \eqref{1} with $(\mu,\nu)=(|\mu_n|,|\nu_n|)$.  (Condition \req{subcr0} implies that every pair of measures is good.)
By comparison principle we have
 $$|u_n|\leq U_n\leq v_n+v'_n\q a.e. $$
 Thus
$$|g\circ u_n|\leq g\circ U_n\leq g\circ(v_n+v'_n)\leq C(|v_n+v'_n|^q+1)\leq C'(|v_n|^q+|v'_n|^q+1).$$

By classical estimates
$$\norm{\BBG(|\mu_n|)}_{L^p(\Gw;\gr)}\leq c_p\norm{\mu_n}_{\CM(\Gw;\gr)} \forevery p\in [1,(N+1)/(N-1))$$
and
$$\norm{\BBP(|\nu_n|)}_{L^p(\Gw;\gr)}\leq c_p\norm{\nu_n}_{\CM(\bdw)} \forevery p\in [1,(N+1)/(N-1)).$$
Hence, $\{v_n\}$ and $\{v'_n\}$ are bounded in $L^p(\Gw;\gr)$ for every $p$ as above. This in turn implies that they are uniformly integrable in each of these spaces. It follows that $\{g\circ u_n\}$ is uniformly integrable in $L^1(\Gw;\gr)$. Since $u_n\to u$ in $L^1(\Gw)$ there exists a subsequence  $\{u_{n_k}\}$ that converges a.e. to $u$.
Therefore $g\circ u_{n_k}\to g\circ u$ .in $L^1(\Gw;\gr)$. As the limit does not depend on the \sseq we conclude that
$g\circ u_n\to g\circ u$ in $L^1(\Gw;\gr)$.

By \req{subcr2} $\{g\circ u_n\}$ is bounded  in $\CM(\Gw;\gr)\}$; therefore a \sseq (still denoted $\{g\circ u_n\}$) converges weakly in $\CM(\bar\Gw)$ to a measure $\gl$.
As $u_n$ is a weak solution of \eqref{26},
\begin{equation}\label{vgf-3}
\int_{\Gw}\big(-u_n\Gd\vgf+(g\circ u_n)\vgf\big)dx=\int_{\Gw}\vgf \  d\mu_n-\int_{\bdw}\frac{\prt\vgf}{\prt\bf n}d\nu_n  \q  \forall\vgf\in C^{2}_{0}(\bar\Gw).
\end{equation}
By \rlemma{prtn} and \req{subcr2},
$$\int_{\Gw}(g\circ u_n)\vgf \ dx\to\int_{\Gw}(g\circ u)\vgf\ dx$$
and
$$\int_{\Gw}\vgf\ d\mu_n\to\int_{\Gw}\vgf\ d\mu_{int}-\int_{\bdw}\frac{\prt\vgf}{\prt\bf n}\ d\tau.$$

Therefore taking the limit in \eqref{vgf-3}, we obtain
$$\int_{\Gw}\displaystyle\left(-u\Gd\vgf+(g\circ u)\vgf\right)dx=\int_{\Gw}\vgf\ d\mu_{int}-\int_{\bdw}\frac{\prt\vgf}{\prt \bf n}d(\nu+\tau\chr{\bdw}).$$
\end{proof}

\section{\bf Negligible measures}

\bth{negligible} Assume that $g\in\CG_0$ is convex and satisfies the $\Gd_2$ condition.
Let $\{(\mu_n,\nu_n)\}$ and $\{(\tl \mu_n,\nu_n)\}$ be $g$-good sequences converging weakly in $\bar\Gw$  to $(\tau,\nu)$ and $(\tl\tau,\nu)$ respectively. Assume that, for every $n\geq 1$, $(|\mu_n|,|\nu_n|)$ and $(|\tl\mu_n|,|\nu_n|)$ are in  $\CM^g(\bar\Gw)$.

Let $u_n$ (resp $\tl u_n$) be the solution of \req{1} with $(\mu,\nu)=(\mu_n,\nu_n)$ (resp. $(\tl\mu_n,\nu_n)$). Assume that
$$u_n\to u,\q \tl u_n\to\tl u \q \text{in }L^1(\Gw)$$
and let $(\mu^*,\nu^*)$ and $(\tl\mu^*,\tl\nu^*)$ denote the reduced limits of $\{(\mu_n,\nu_n)\}$ and $\{(\tl\mu_n,\nu_n)\}$ respectively.

Assume that a \sseq of
$\{\gr|\tl\mu_{n}-\mu_{n}|\}$  converges weakly in $\CM(\bar\Gw)$ to a measure $\Gl$ \sth $\Gl\chr{\bdw}$ is negligible.
\bcom +++++++++
Assume that,
\begin{equation}\label{neg1'}
    \Gl\chr{\bdw}\q\text{ is negligible.}
\end{equation}
Under these assumptions, if
\begin{equation}\label{neg1}
 \mu_n\leq\tilde\mu_n,\q (\tl\tau-\tau)\chr{\bdw} \;\text{ is negligible}
\end{equation}
++++++++++\end{comment}
Then
\begin{equation}\label{neg2}
   \nu^{*}=\tl\nu^{*}.
\end{equation}
\es

\begin{proof} First we prove the result in the case that
\begin{equation}\label{neg1}
   \mu_n\leq\tilde\mu_n.
\end{equation}

This condition implies that $u_n\leq \tl u_n$ and \consy $\nu^*\leq \tl\nu^*$.
By Lemma \ref{l:comp2},
\begin{equation}\label{comp-nu}
0\leq\tl\nu^{*}-\nu^{*}\leq (\tl\tau-\tau)\chr{\bdw}.
\end{equation}

Observe that
\begin{equation}\label{abs-nu}
 |\nu^*|\in \CM^g(\bdw),\q  |\tl\nu^*|\in \CM^g(\bdw),\q \tl\nu^*-\nu^*\in \CM^g(\bdw).
\end{equation}
Passing to a \sseq we may assume that:
(a)$\{(|\mu_n|,|\nu_n|)\}$ possesses a reduced limit $(\bar\mu,\bar\nu)$ and (b) $\{\gr|\tl\mu_{n}-\mu_{n}|\}$  converges weakly in $\CM(\bar\Gw)$ to a measure $\Gl$  \sth $\Gl\chr{\bdw}$ is negligible.

 Thus $(\bar\mu,\bar\nu)\in \CM^g(\bar\Gw)$; since both measures are positive it follows that $\bar\mu\in \CM^g(\Gw)$
and $\bar\nu\in \CM^g(\bdw)$. Clearly $|\nu^*|\leq \bar\nu$; therefore $|\nu^*|\in \CM^g(\bdw)$. Similarly $|\tl\nu^*|\in \CM^g(\bdw)$. In view of our assumptions on $g$, these facts imply that $\tl\nu^*-\nu^*\in \CM^g(\bdw)$.

Since $ (\tl\tau-\tau)\chr{\bdw}$ is negligible while $\tl\nu^{*}-\nu^{*}$ is a non-negative measure in $\CM^g(\bdw)$, \req{comp-nu} implies that $\nu^*=\tl \nu^*$.

Next we drop condition \req{neg1}. Without loss of generality we may assume that the  entire \seq $\{\gr|\tl\mu_{n}-\mu_{n}|\}$  converges weakly in $\CM(\bar\Gw)$ to $\Gl$.

Put $\gg_n:=\mu_{n}+|\tl\mu_{n}-\mu_{n}|$.
Since $g$  is super additive (as a consequence of the convexity assumption and the fact that $g(x,0)=0$) and satisfies the $\Gd_2$ condition $|\gg_n|\in \CM^g(\Gw)$. Since $|\nu_{n}|$ is a good measure it follows that $(|\gg_n|,|\nu_{n}|)\in \CM^g(\bar\Gw)$. Passing to a \sseq we may assume that$\{(\gg_n,\nu_{n})\}$ converges weakly in $\bar\Gw$ and possesses a reduced limit $(\gg^*,\nu^*_1)$.

Note that
$$\mu_{n}\leq \gg_n,\q \tl\mu_{n}\leq \gg_n \forevery n\geq 1.$$
Furthermore,
$$\gr(\gg_n-\mu_{n})\rightharpoonup\Gl.$$
Therefore, by the first part of the proof, applied  to the sequences $\{(\gg_n,\nu_{n})\}$ and $\{(\mu_{n},\nu_{n})\}$ we obtain,
$$\nu^*=\nu^*_1.$$

Next observe that $$|\gg_n-\tl\mu_{n}|\leq 2|\tl\mu_{n}-\mu_{n}|.$$
Consider a \sseq of $\{\gr|\gg_n-\tl\mu_{n}|\}$ that converges weakly in $\CM(\bar\Gw)$ to a measure $\Gl'$. Then $\Gl'\leq 2\Gl$ and, as $\Gl\chr{\bdw}$ is negligible, it follows that
$\Gl'\chr{\bdw}$ is negligible.
Applying the first part of the proof to the sequences  $\{(\gg_n,\nu_{n})\}$ and $\{(\tl\mu_{n},\nu_{n})\}$ we obtain,
$$\tl\nu^*=\nu^*_1.$$
Combining these facts we obtain \req{neg2}.
\end{proof}
\vskip 2mm

\Remark If all the measures are non-negative and $\mu_n\leq \tl\mu_n$ then the conclusion of the theorem is valid for every $g\in \CG_0$, i.e., convexity and the $\Gd_2$ condition are not needed. Indeed in this case
$\nu^*$ and $\tl\nu^*$ are non-negative and $\nu^*\leq \tl \nu^*$. Furthermore, by definition, the reduced limits belong  to $\CM^g(\bar\Gw)$. As the measures are non-negative this implies that $\nu^*$ and $\tl \nu^*$ are in $\CM^g(\bdw)$. These facts imply that $\tl\nu^*-\nu^*$ is a non-negative good measure. As $\tl\tau-\tau$ is negligible, \req{comp-nu} implies that $\nu^*=\tl\nu^*$.

\bcom @@@@@@@@@@@@@@@@@@@@@@@@@@@@@@@@@@@@@@@@@@@@@
Next we claim that if either (a) or (b) holds,
\begin{equation}\label{neg3}
\tl\nu^*-\nu^*\in \CM^g(\bdw).
\end{equation}
Let $v$ (resp. $\tl v$) be the weak solution of \req{1} with data $(0,\nu^*)$ (resp. $(0,\tl\nu^*)$. Then $v\leq \tl v$.
If condition (b) holds then (in view of the superadditivity of $g$),
$$g\circ(\tl v-v)\leq (g\circ \tl v)-(g\circ v).$$
It follows that $g\circ(\tl v-v)\in L^1(\Gw,\rho)$ and that $\tl v- v$ is a subsolution of \req{1} with data $(0,\tl\nu^*-\nu^*)$.  On the other hand, as $\nu^*\geq 0$, $\tl v$ is a supersolution of this problem. Therefore, by  \cite[Theorem 2.2.4]{MV} (or \cite{Mon,P}), this problem has a solution, i.e., $\tl\nu^{*}-\nu^{*}\in\CM^{g}(\bdw)$.

Next assume that condition (a) holds. Then, by superadditivity of $g$
$$ 0\leq g\circ(\tl v-v)\leq g\circ (\tl v^++v^-)- g\circ (\tl v^-+v^+).$$
If $g$ satisfies the $\Gd_2$ condition then both terms on the right hand side belong to $L^1(\Gw;\gr)$. Thus $\tl v-v$ is a subsolution of \req{1} with data $(0,\tl\nu^*-\nu^*)$.
We have shown above that $|\nu^*|\in \CM^g(\bdw)$ and $|\tl\nu^*|\in \CM^g(\bdw)$. Therefore problem

{\bf Claim:} $\tl\nu^{*}-\nu^{*}\in\CM^{g}(\bdw)$.\\
By Proposition \ref{p:class}(i) it follows that $\tl\nu^{*}, \nu^{*}\in\CM^g(\bdw)$. Let $v$ is a weak solution of \eqref{1} with
$(\mu,\nu)=(0,\nu^{*})$ and similarly $w$ is a weak solution of \eqref{1} with $(\mu,\nu)=(0,\tl\nu^{*})$. By \eqref{comp-nu}, $0\leq v\leq w$. Therefore $g\circ(w-v)\leq g\circ w-g\circ v$ and $g\circ(w-v)\in L^1(\Gw,\rho)$. Thus $w-v$ (resp $w$) is a subsolution (resp supersolution) of \eqref{1} with $(\mu,\nu)=(0, \tl\nu^{*}-\nu^{*})$. Hence by \cite[Theorem 2.2.4]{MV}, $\tl\nu^{*}-\nu^{*}\in\CM^{g}(\bdw)$.

@@@@@@@@@@@@@@@@@@
(ii)Assume $\mu_n, \nu_n\leq 0$ and $-\tau\chr{\bdw}$ is negligible. Now as $g(x,\cdot)$ is a odd function, we can easily check that $-u_n$ is a weak solution of \eqref{26} where $(\mu_n,\nu_n)$ is replaced by $(-\mu_n,-\nu_n)$. Similarly $-v_n$ is a weak solution of \eqref{bvpn6} where $\mu_n$ is replaced by $-\mu_n$ and $-w_n$ is a weak solution of \eqref{bvpn7} where $\nu_n$ is replaced by $-\nu_n$. Now proceeding the same way as in part (i), we can prove that trace $u\geq$trace $v+$trace $w$. On the other hand since $(-\tau\chr{\bdw})$ is negligible, by the same analysis of part (i) yields trace($-v)=0$ which implies\\ trace $v=0$. Therefore trace $u\geq$trace $w$. As $\mu_n\leq 0$, we also have $u_n\leq w_n$. Hence trace $u=$trace $w=\nu^{\#}$.

\bcor{neg} Let $g\in \CG_0$ and assume that $g$ is convex and satisfies the $\Gd_2$ condition. Assume all the assumptions of the theorem  except for \req{neg1} which is replaced by condition \req{neg1'} stated below.

There exists a \sseq, say
$\{\gr|\tl\mu_{n_k}-\mu_{n_k}|\}$, that converges weakly in $\CM(\bar\Gw)$ to a measure $\Gl$. Assume that,
\begin{equation}\label{neg1'}
    \Gl\chr{\bdw}\q\text{ is negligible.}
\end{equation}
Then \req{neg2} holds.
\es
\begin{proof} Put $\gg_k=\mu_{n_k}+|\tl\mu_{n_k}-\mu_{n_k}|$.
Since $g$ is superadditive and satisfies the $\Gd_2$ condition $|\gg_k|\in \CM^g(\Gw;\gr)$. Since $|\nu_{n_k}|$ is a good measure it follows that $(|\gg_k|,|\nu_{n_k}|)\in \CM^g(\bar\Gw)$. Passing to a \sseq we may assume that$\{(\gg_k,\nu_{n_k})\}$ converges weakly in $\Gw$ and possesses a reduced limit $(\gg^*,\nu^*_0)$.

Note that
$$\mu_{n_k}\leq \gg_k,\q \tl\mu_{n_k}\leq \gg_k \forevery k\geq 1.$$
Furthermore,
$$\gr(\gg_k-\mu_{n_k})\rightharpoonup\Gl.$$
Therefore applying the theorem to the sequences $\{(\gg_k,\nu_{n_k})\}$ and $\{(\mu_{n_k},\nu_{n_k})\}$ we obtain,
$$\nu^*=\nu^*_0.$$

Next observe that $$|\gg_k-\tl\mu_{n_k}|\leq 2|\tl\mu_{n_k}-\mu_{n_k}|.$$
Consider a \sseq of $\{\gr|\gg_k-\tl\mu_{n_k}|\}$ that converges weakly in $\CM(\bar\Gw)$ to a measure $\Gl'$. Then $\Gl'\leq 2\Gl$ and \req{neg1'} implies  that
$$\Gl'\chr{\bdw}\text{ is negligible.}$$
Therefore, passing again to a \sseq, we may apply the theorem to $\{(\gg_k,\nu_{n_k})\}$ and $\{(\tl\mu_{n_k},\nu_{n_k})\}$. It follows that
$$\tl\nu^*=\nu^*_0.$$
Combining these facts we obtain \req{neg2}.
\end{proof}
@@@@@@@@@@@@@@@@@@@@@@@@@@@\end{comment}

\end{document}